\pdfoutput=1

\documentclass[twocolumn]{autart}    

\usepackage{graphicx} 
\usepackage{tikz}

\usepackage{bm}
\usepackage{enumerate}
\usepackage{commath}
\usepackage{graphicx}
\graphicspath{{Pictures/}}
\usepackage{amsmath}
\usepackage{breqn}
\usepackage{cite}
\usepackage{booktabs}
\usepackage{empheq}
\usepackage{amsfonts}
\usepackage{amssymb}

\usepackage{xcolor}
\usepackage{mathrsfs}
\usepackage{multirow}
\usepackage{changes}
\usepackage{color}

\usepackage{accents}
\usepackage{thmtools}
\usepackage{thm-restate}
\usepackage{float}
\usepackage{comment}
\usepackage[normalem]{ulem}
\usepackage{algorithm}
\usepackage{algpseudocode}
\usepackage[normalem]{ulem}
\newtheorem{lemma}{Lemma}
\newtheorem{definition}{Definition}
\newtheorem{assumption}{Assumption}
\newtheorem{proposition}{Proposition}
\newtheorem{remark}{Remark}
\newcommand*{\QEDA}{\hfill\ensuremath{\blacksquare}}
\usepackage{algcompatible}
\usepackage{cancel}

\begin{document}

\begin{frontmatter}

\title{Achieving violation-free distributed optimization under coupling constraints\thanksref{footnoteinfo}} 

\thanks[footnoteinfo]{This work was supported by the Swedish Research Council, the Knut and Alice Wallenberg Foundation, the Swedish Foundation for Strategic Research, and an NSERC Postdoctoral Fellowship.}

\author[KTH]{Changxin Liu}\ead{changxin@kth.se},    
\author[KTH]{Xiao Tan}\ead{xiaotan@kth.se},               
\author[KTH]{Xuyang Wu}\ead{xuyangw@kth.se},  
\author[KTH]{Dimos V. Dimarogonas}\ead{dimos@kth.se},  
\author[KTH]{Karl H. Johansson}\ead{kallej@kth.se}  

\address[KTH]{School of Electrical Engineering and Computer Science, KTH Royal Institute of Technology, and Digital Futures, Stockholm, Sweden}  


\begin{abstract}                          
Constraint satisfaction is a critical component in a wide range of engineering applications, including but not limited to safe multi-agent control and economic dispatch in power systems. This study explores violation-free distributed optimization techniques for problems characterized by separable objective functions and coupling constraints.
First, we incorporate auxiliary decision variables together with a network-dependent linear mapping to each coupling constraint. 
For the reformulated problem, we show that the projection of its feasible set onto the space of primal variables is identical to that of the original problem, which is the key to achieving all-time constraint satisfaction.
Upon treating the reformulated problem as a min-min optimization problem with respect to auxiliary and primal variables, we demonstrate that the gradients in the outer minimization problem have a locally computable closed-form.
Then, two violation-free distributed optimization algorithms are developed and their convergence under reasonable assumptions {is} analyzed. Finally, the proposed algorithm is applied to implement a control barrier function based controller in a distributed manner, and the results verify its effectiveness.
\end{abstract}

\end{frontmatter}

\section{Introduction}
Distributed optimization problems over networks, such as economic dispatch in power networks and coordination in multi-agent systems, typically involve decision variables that are subject to coupling constraints. {
These constraints typically represent the shared resources of the parties involved and cannot be violated because of physical limitations.
Due to its paramount importance, the design of distributed optimization algorithms for constraint-coupled problems has received increasing attention lately (see, for example,  \cite{notarnicola2019constraint,li2020distributed,tan2021distributed,wu2022distributed,camisa2021distributed}).

An essential requirement in designing distributed optimization algorithms for such problems is all-time satisfaction of the coupling constraints. This is because the constraint satisfaction is critical to ensuring that the system operates safely and effectively \cite{ames2016control}. Furthermore, all-time satisfaction guarantees that the iterative algorithm ends up with a feasible and safe-to-implement solution whenever it is stopped. Nevertheless, most existing results only have asymptotic feasibility guarantee, that is, the constraint violation vanishes asymptotically \cite{notarnicola2019constraint,li2020distributed,camisa2021distributed}. Projecting such solutions onto the feasible set in presence of coupling constraints is usually not permissible in a distributed scenario. To this end, this work aims at developing violation-free distributed algorithms for optimization problems with separable objective functions and coupling constraints.

}

\subsection{Related works}
{
The most well-known methodology for solving constraint-coupled distributed optimization problems is dual decomposition \cite{boyd2011distributed}. 
For instance, the authors in \cite{falsone2017dual,nedic2009approximate,liu2020unitary} considered the Lagrangian dual of the original problem and developed algorithms for distributed optimization problems with general separable constraints.
Based on a similar saddle-point formulation, the authors in \cite{mateos2016distributed} proposed distributed primal-dual algorithms, where the primal variables are updated via one-step gradient descent rather than solving subproblems.
{Notably, the authors in \cite{notarnicola2019constraint} considered a relaxed version of the constraint-coupled optimization problem, based on which a distributed algorithm was designed to ensure the last iterate convergence of variables without averaging steps.
Utilizing the same methodology, the authors in \cite{wang2023distributed} developed an improved optimization algorithm for strongly convex problems.} 
For sparsely coupled equality constraints, the authors in \cite{alghunaim2019proximal} introduced a structure-aware Lagrangian dual formulation, based on which the overall algorithm can be made more efficient.
However, the primal iterates by all the above algorithms only satisfy the coupling constraints asymptotically, that is, the constraint violation asymptotically vanishes. Therefore, they may not generate feasible solutions within finite time, and are not implementable to safety-critical systems.}


Another strategy for designing constraint-coupled distributed optimization is primal decomposition \cite{wu2021new,tan2021distributed,wu2022distributed}, which enforces all-time constraint satisfaction. Among these methods, \cite{tan2021distributed,wu2021new} are originated from the right-hand side allocation strategy in \cite{bertsekas2016nonlinear}, and involve  decomposing the coupling constraint with the help of additional variables and updating the additional variables along a direction such that the constraint stays feasible. Particularly, the authors in \cite{tan2021distributed} studied problems with $1$-dimensional coupling constraint, and presented a continuous-time distributed optimization algorithm based on finite-time consensus-seeking protocols. For problems subject to multi-dimensional coupling constraints, the authors in \cite{wu2021new,wu2022distributed} developed distributed feasible methods, where the update direction simultaneously decreases the objective function value and keeps the coupling constraint feasible. However, such update direction is found by each agent upon directly exchanging objective \cite{wu2021new} or gradient information \cite{wu2022distributed}, which poses a potential risk of privacy leakage \cite{zhu2019deep}. It is noteworthy that only convergence to a neighborhood of the optimal solution is guaranteed in \cite{wu2022distributed,wu2021new}. 

{
Some other notable methods achieving all-time constraint satisfaction include \cite{turan2022safe,doostmohammadian2022distributed}. The authors in \cite{turan2022safe} considered a special class of constraint-coupled optimization problems, where the coefficients of the linear constraints belong to $\{0, 1\}$. A safety margin is constructed to avoid constraint violation by tightening the original constraint properly. In \cite{doostmohammadian2022distributed}, $1$-dimensional distributed optimization with equality constraint is studied. For this specific problem, the authors exploited the property that the local gradients evaluated over optimal solution are identical, and developed a distributed algorithm that averages the local gradients iteratively. The recursive satisfaction of coupling constraint follows from the property of Laplacian used for averaging. 

To summarize, most distributed optimization algorithms
handling coupling constraints in the literature cannot provide all-time constraint satisfaction, and those that can provide either focus on problems with specific structure such as one-dimensional variables and cannot be easily extended to handle a more general class of problems, or suffer from inexact convergence.


}





\subsection{Our contribution}
In this work, we consider networked optimization problems with separable convex objective functions and coupling multi-dimensional constraints in the form of both equalities and inequalities. Compared to existing constraint-coupled distributed optimization \cite{wu2022distributed,notarnicola2019constraint,falsone2017dual}, the proposed algorithms possess the following notable characteristics:
\begin{itemize}
    \item[i)] They produce violation-free solutions whenever they are terminated, while also converging to precise solutions with an explicit rate guarantee.
    \item[ii)] They leverage the inherent structure of the coupling constraints, leading to enhanced communication efficiency and convergence performance. 
\end{itemize}


To accomplish these desirable objectives, we reformulate multiple constraints by introducing auxiliary variables with a particular network-dependent linear transformation. This reformulation enables the decomposition of the problem, making it amenable to distributed solutions.
The auxiliary variables corresponding to a certain constraint are introduced to an agent only if it is affected by the constraint, leading to a sparse and efficient formulation.

{
Subsequently, the reformulated problem is approached as a min-min optimization scenario, where the auxiliary and primal variables are optimized separately, and examined through sensitivity analysis.
In particular, we show that the gradients of the objective function in the outer minimization are network-dependent affine transformations of Karush-Kuhn-Tucker (KKT) multipliers of the inner problem under mild conditions, and can be locally computed by agents.
}
Provided that the local objective is strongly convex, we quantify the Lipschitz constants of the gradients, which facilitates the use of the accelerated dual averaging algorithm \cite{cohen2018acceleration} in solving the reformulated problem. For general convex objectives, additional coordinate constraints are imposed on the auxiliary variables, which ensures the boundedness of the gradients. Based on this, the reformulated problem is solved by gradient descent with convergence guarantee.

Finally, the proposed algorithm is tested on a constrained consensus-seeking system under a control barrier function (CBF) based controller \cite{xu2015robustness,ames2016control}, where a quadratic programming problem with sparse coupling constraints is solved at each sampling time. The results verify the effectiveness of the algorithm. 

\subsection{Paper organization and notation}
The structure of this work is outlined as:
Section \ref{sec:2} introduces the networked optimization problem with coupling constraints. In Section \ref{sec:3}, we present a reformulation of the problem and propose a violation-free distributed optimization algorithm for strongly-convex problems, and provide convergence rate results. We extend the algorithm and analysis to general convex problems in Section \ref{sec:4}. Section \ref{sec:5} presents the results of numerical experiments, and finally, Section \ref{sec:6} concludes this work.

$\mathbb{R}$, $\mathbb{R}^d$, and $\mathbb{R}^{n\times d}$ represent the 1-dimensional, d-dimensional, and $n\times d$ Euclidean space, respectively. $\mathbb{N}$ is the set of natural numbers.
$\text{blkdiag}(g_1, g_2, \dots, g_n )$ denotes a block diagonal matrix with its diagonal blocks $g_1, g_2, \dots, g_n$, where $g_i, i = 1,\dots,n$ can be either a vector or a matrix, and $\text{col}(h_1, h_2, \dots, h_n) = \begin{bmatrix}
   h_1^T, h_2^T, \dots, h_n^T
\end{bmatrix}^T$, where $h_i, i = 1,\dots,n$ is a column vector. We denote by $\text{mat}(A_1,\dots,A_n)=\begin{bmatrix}
    A_1& \cdots & A_n
\end{bmatrix}$, where $A_i, i = 1,\dots,n$ can be either a column vector or a matrix. Given two sets $\mathcal{A}$ and $\mathcal{B}$, $\oplus$ stands for the Minkowski sum defined by
$
    \mathcal{A} \oplus \mathcal{B}:=\{a+b|a\in \mathcal{A}, b\in\mathcal{B}  \}.
$

\section{Problem statement}
\label{sec:2}
\subsection{Basic setup}
Consider a class of linearly constrained optimization problems given by
\begin{equation} \label{eq:centralized_problem}
    \begin{split}
        &\min_{x_1,\dots,x_N} \sum_{i=1}^N f_i(x_i) \\
        & \text{s.t.} \quad \sum_{i =1}^N A_i x_i +b_{i} \leq 0, \\
        & \quad \quad \,\,\,  \sum_{i =1}^N E_i x_i +g_{i} = 0
    \end{split}
\end{equation}
where
$x_i\in\mathbb{R}^{d_i}$, $d_i\in\mathbb{N}$, $A_{i} \in \mathbb{R}^{ M\times d_i}$, ${ b}_{i} \in \mathbb{R}^M$, $E_{i} \in \mathbb{R}^{Q\times d_i}$,
${g}_{i} \in \mathbb{R}^Q$, and $f_i:\mathbb{R}^{d_i}\rightarrow \mathbb{R}$ are local to each agent $i=1,\dots N$. We denote an optimal solution to \eqref{eq:centralized_problem} by $x_i^*, i=1,\dots,N$. 



Each pair $(A_i, b_i)$ (resp. $(E_i, g_i)$) has $M$ (resp. $Q$) rows.
Denote by $A_i^{[m]}$ the $m$-th row of $A_i$ and by $b_i^{[m]}$ the $m$-th entry of $b_i$. Similarly, let $E_i^{[q]}$ and $g_i^{[q]}$ be the $q$-th row and coordinate of $E_i$ and $g_i$, respectively. Some of the rows can take zero values in the case of sparse coupling constraints. 
{Hereafter, we denote by $\mathcal{V}^{[m]}=\{j: A_j^{[m]}\ne 0 \,\, \mathrm{or} \,\, b_j^{[m]}\ne 0 \}$ (resp. $\mathcal{V}^{[M+q]}=\{j: E_j^{[q]}\ne 0 \,\, \mathrm{or} \,\, g_j^{[q]}\ne 0 \}$) the set of agents that are influenced by the $m$-th inequality constraint (resp. $q$-th equality constraint), and by $\mathcal{I}_i$ (resp. $\mathcal{E}_i$) the set of indexes of inequality (resp. equality) constraints that influence agent $i$'s decision variable.} 



\subsection{Communication network}
The communication among the agents is described by an undirected graph $\mathcal{G}=(\mathcal{V},\mathcal{S})$, where $\mathcal{V}$ is the set of agent indexes and $\mathcal{S}$ is the set of unordered pairs $(i,j),i,j\in\mathcal{V}$. 
Each pair $(i,j) \in \mathcal{S}$ represents a communication channel between $i$ and $j$, and $\mathcal{S}\subseteq\mathcal{V}\times\mathcal{V}$ represents the set of communication channels.
We denote by $\mathcal{N}_{i}$ the cluster of $i$'s neighbors including itself, i.e., the set of agents that $i$ can communicate with, {i.e., $\mathcal{N}_i=\{j\in\mathcal{V}: (i,j)\in\mathcal{S}\}\cup\{i\}$. 
}

For each scalar constraint indicated by the corresponding rows of $A_i$'s and $E_i$'s, the involved agents as well as the links between them in the communication $\mathcal{G}$ form an induced subgraph \cite{mesbahi2010graph}. We denote by $\mathcal{G}^{[m]}=(\mathcal{V}^{[m]},\mathcal{S}^{[m]})$ (resp. $\mathcal{G}^{[M+q]}=(\mathcal{V}^{[M+q]},\mathcal{S}^{[M+q]})$) the graph induced by the $m$-th inequality (resp. $q$-th equality) constraint. {Specifically, $\mathcal{S}^{[m]}=\{(i,j): (i,j)\in\mathcal{S}, i,j\in \mathcal{V}^{[m]}\}$ and the definition of $\mathcal{S}^{[M+q]}$ is similar.} Accordingly, for each $\mathcal{G}^{[l]}, l=1,\dots, M+Q$, we denote by $\mathcal{N}_{i}^{[l]}$ the cluster of $i$'s neighbors {including itself}. The following example is presented to illustrate the communication network.

{
\emph{Example}. Consider problem \eqref{eq:centralized_problem} with $N=4$, $d_i=3, \forall i$, $M=Q=1$,
\begin{equation*}
    \begin{pmatrix}
        A_1\\
        A_2\\
        A_3\\
        A_4
    \end{pmatrix} = \begin{pmatrix}
        1 & 0 & 0\\
        0 & 0 & 0\\
        0 & 0 & 0\\
        1 & 1 & 1
    \end{pmatrix},\qquad     \begin{pmatrix}
        E_1\\
        E_2\\
        E_3\\
        E_4
    \end{pmatrix} = \begin{pmatrix}
        1 & 0 & 0\\
        1 & 1 & 0\\
        0 & 1 & 0\\
        0 & 0 & 0
    \end{pmatrix},
\end{equation*}
and the communication graph in Fig. \ref{fig:graph}.
\begin{figure}[!htb]
    \centering
    \begin{tikzpicture}[node distance=2cm, auto]
        \node[circle, draw] (1) {1};
        \node[circle, draw, right of=1] (2) {2};
        \node[circle, draw, below of=2] (3) {3};
        \node[circle, draw, left of=3] (4) {4};
        \draw[-] (1) -- (2);
        \draw[-] (2) -- (3);
        \draw[-] (3) -- (4);
        \draw[-] (4) -- (1);
        \draw[-] (1) -- (3);
    \end{tikzpicture}
    \caption{A graph with 4 nodes.}
    \label{fig:graph}
\end{figure}
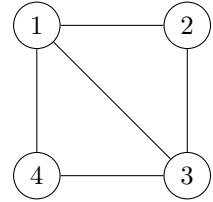
Then,
\begin{align*}
    &\mathcal{V}^{[1]} = \{1,4\},~~\mathcal{S}^{[1]}=\{(1,4)\},\\
    &\mathcal{V}^{[2]} = \{1,2,3\}, ~~\mathcal{S}^{[2]}=\{(1,2), (2,3), (1,3)\}.
\end{align*}

The following assumption is made for the communication network.


\begin{assumption}
\label{assump:network}
For each $l=1,\dots, M+Q$,
    $\mathcal{G}^{[l]}$ is undirected and connected.
\end{assumption}

Assumption \ref{assump:network} is not restrictive in the sense that if all $A_i^{[m]},m=1,\dots,M$ and $E_i^{[M+q]},q=1,\dots,Q$ are non-zero then Assumption \ref{assump:network} holds when $\mathcal{G}^{[l]} = \mathcal{G},l=1,\dots, M+Q$ is connected.

We associate each subgraph $\mathcal{G}^{[l]}$ with a \textit{symmetric}, nonnegative weight matrix $P^{[l]}=[p^{[l]}_{ij}]\in\mathbb{R}^{\lvert \mathcal{V}^{[l]} \rvert \times \lvert \mathcal{V}^{[l]} \rvert }$, which will be used for the agents to weigh the exchanged information from their neighbors. The weights $p^{[l]}_{ij}$ satisfy the following conditions
\begin{equation}\label{eq:row_stochastic}
    p_{ij}^{[l]} = \left\{\begin{aligned}
    >&0,  & j\in\mathcal{N}_{i}^{[l]} \\
    &0, & \text{otherwise}
    \end{aligned}\right., \quad \sum_{j=1}^N p_{ij}^{[l]}=1\,\,\, \forall i. 
\end{equation}
Eq. \eqref{eq:row_stochastic} indicates the double stochasticity of $P^{[l]}$, and this requirement can be fulfilled by several well-known protocols, e.g., Metropolis-Hastings rule \cite{boyd2004fastest}. We note that if $\mathcal{G}^{[l]}$ is connected, there holds $\text{Null}(I-P^{[l]})=\text{Span}({\bf 1})$ and $\text{Range}(I-P^{[l]}) =  \{ {z}\in \mathbb{R}^{\lvert \mathcal{V}^{[l]} \rvert}: {\bf 1}^T {z} = 0\}$, where $\text{Null}(\cdot)$ and $\text{Range}(\cdot)$ denotes the null space and the range of a linear map \cite{shi2015extra}.

{
\begin{remark}
    Assumption \ref{assump:network} imposes a connectivity condition on each subgraph $\mathcal{G}^{[l]}$, such that only the agents affected by the $l$-th constraint will participate in allocating this constraint by communication and local computation. Therefore, the computation and communication of the overall algorithm can be made more efficient. 
\end{remark}}





\section{Problem reformulation and algorithm}
\label{sec:3}
This section presents a novel reformulation of the optimization problem in \eqref{eq:centralized_problem}. Based on sensitivity analysis, we show that the gradients of the objective function of the reformulated problem can be locally computed. Upon exploiting this feature, we develop an accelerated distributed optimization algorithm for solving \eqref{eq:centralized_problem}.

\subsection{Problem reformulation}


Before proceeding to the reformulation, we introduce the following notation
\begin{equation*}
\begin{split}
    &{\bf A}^{[m]} = \text{blkdiag}(\{A_{l}^{[m]} \}_{l \in \mathcal{V}^{[m]}}),\,\, { b}^{[m]} = \text{col}(\{b_{l}^{[m]} \}_{l \in \mathcal{V}^{[m]}})  \\
    & {\bf E}^{[q]} = \text{blkdiag}(\{E_{l}^{[q]} \}_{l \in \mathcal{V}^{[M+q]}}),  \,\, { g}^{[q]} = \text{col}(\{g_{l}^{[q]} \}_{l \in \mathcal{V}^{[M+q]}}) \\
    &x^{[m]}= \text{col}(\{x_l\}_{l\in\mathcal{V}^{[m]}}), \quad \quad  \,\,\, x^{[M+q]}= \text{col}(\{x_l\}_{l\in\mathcal{V}^{[M+q]}}).
    \end{split}
\end{equation*}
We incorporate a vector of slack variables to each constraint in \eqref{eq:centralized_problem}, denoted by $y^{[l]}=\text{col}(\{ y^{[l]}_i \}_{i \in \mathcal{V}^{[l]}}) \in \mathbb{R}^{\lvert \mathcal{V}^{[l]} \rvert}, l=1,\dots, M+Q$. Further denoting ${y}=\text{col}(\{y^{[m]}\}_{m=1,\dots,M+Q})$ and ${ x}=\text{col}(\{x_i\}_{i=1,\dots,N})$,
we arrive at the following new problem
\begin{equation}\label{transformed_P}
\begin{split}
&\min_{x,y } \sum_{i=1}^N f_i(x_i) \\
    & \text{s.t.} \quad 
   { \bf A}^{[m]}  {  x}^{[m]} + \left(I- P^{[m]}\right) { y}^{[m]} + { b}^{[m]}\leq 0, \\
  & \quad  \,\,\quad  \forall m=1,\dots M \\
 & \quad \quad \,\,  {\bf E}^{[q]} {  x}^{[M+q]}  + \left(I- P^{[M+q]}\right) { y}^{[M+q]} + { g}^{[q]}= 0,\\
  & \quad  \,\,\quad  \forall q=1,\dots Q
\end{split}
\end{equation}
where $P^{[m]}$ is the weight matrix defined in \eqref{eq:row_stochastic}. {The utilization of the linear mapping $(I-P^{[m]})$ for slack variables $y^{[m]}$ is motivated by two key factors. Firstly, it facilitates a network-aware decomposition of the problem, while ensuring the preservation of the solution for the original optimization problem \eqref{eq:centralized_problem}, as outlined in Proposition \ref{prop:equivalence}. Secondly, compared to alternative linear mappings such as the graph Laplacian, $(I-P^{[m]})$ possesses a lower matrix norm. This characteristic is advantageous as it allows for the amplification of the step-size and therefore speeds up the convergence process.}

\begin{proposition}\label{prop:equivalence}
    Suppose Assumption \ref{assump:network} holds. The problems in \eqref{eq:centralized_problem} and \eqref{transformed_P} are equivalent in the sense that
    \begin{itemize}
        \item[i)] they share the same objective function,
        \item[ii)] for any feasible solution $({x},{y})$ to \eqref{transformed_P}, $x$ is a feasible solution to \eqref{eq:centralized_problem}, 
        \item[iii)] for any feasible solution $x$ to \eqref{eq:centralized_problem}, there exist some $y$ such that $({x},{y})$ is feasible to \eqref{transformed_P}. 
    \end{itemize}
\end{proposition}
\begin{pf}
        The statement in i) holds trivially.
        
        ii) Given any feasible solution $({x},{y})$ to \eqref{transformed_P}, it holds, for any $m$ that
        \begin{equation*}
        \begin{split}
        &{\bf 1}^T  \left( { \bf A}^{[m]}  {  x}^{[m]} + \left(I- P^{[m]}\right) { y}^{[m]} + { b}^{[m]}\right)  \\
        &= \sum_{i\in\mathcal{V}} {{A}_i^{[m] }}{x}_i  + b_i^{[m]}  \leq 0  
        \end{split}
        \end{equation*}
        and any $q$ that
        \begin{equation*}
                \begin{split}
        &{\bf 1}^T  \left( { \bf E}^{[q]}  {  x}^{[M+q]} + \left(I- P^{[M+q]}\right) { y}^{[M+q]} + { g}^{[q]}\right) \\
        &= \sum_{i\in\mathcal{V}} {{E}_i^{[q] }}{x}_i  + g_i^{[q]}  = 0
                \end{split}
        \end{equation*}
        because of column stochasticity of $P^{[l]}, l=1,\dots, M+Q$.
        Thus, ${x}$ also satisfies the constraint in \eqref{eq:centralized_problem}.  
    
    iii) Given any feasible solution ${x}$ to \eqref{eq:centralized_problem}, we define 
    \begin{equation*}
    \begin{split}
        {r}^{[m]} &= { \bf A}^{[m]}  {  x}^{[m]} +  { b}^{[m]}, \forall m \\
        {v}^{[q]} &= { \bf E}^{[q]}  {  x}^{[M+q]} +  { g}^{[q]}, \forall q.
    \end{split}    
    \end{equation*}
    Then, it holds that $\overline{r}^{[m]} = \frac{{\bf 1}^T {r}^{[m]}}{\lvert \mathcal{V}^{[m]}\rvert}  \leq 0$ and $\overline{v}^{[q]} = \frac{{\bf 1}^T {v}^{[q]}}{\lvert \mathcal{V}^{[M+q]}\rvert}   =0$.
    Note that  $\{ (I-P^{[l]}){y}^{[l]}: {y}^{[l]}\in \mathbb{R}^{\lvert \mathcal{V}^{[l]} \rvert}\} =  \text{Range}(I-P^{[l]}) = \{ {z}\in \mathbb{R}^{\lvert \mathcal{V}^{[l]} \rvert}: {\bf 1}^T {z} = 0\}$.
    Thus, there exists a ${y}^{[l]}$ such that $(I-P^{[l]}){y}^{[l]} + {r}^{[l]} = {\bf 1}\overline{r}^{[l]}$ for $l=1,\dots,M$ and that $(I-P^{[l]}){y}^{[l]} + {v}^{[l]} = {\bf 1}\overline{v}^{[l]}$ for $l=M+1,\dots,M+Q$.
    Based on this, a feasible solution $({x},{y})$ can be constructed for \eqref{transformed_P}. \QEDA
\end{pf}





Given fixed slack variables $\{{ y}^{[l]}\}_{l=1,\dots, M+Q}$, \eqref{transformed_P} can be partitioned and assigned to each agent. 
In particular, the local optimization problem for agent $i$ is
\begin{equation}\label{transformed_P_i}
\begin{split}
&\min_{x_i}  f_i(x_i) \\
    & \text{s.t.} \,\,
    A_{i}^{[m]} { x}_i   +  y^{[m]}_{i} - \sum_{j=1}^N p_{ij}^{[m]}y^{[m]}_{j}  + { b}^{[m]}_{i}\leq 0, \, \forall m\in\mathcal{I}_i \\
    & \quad \, \, \,   E_i^{[q]} { x}_i   + y^{[M+q]}_{i} - \sum_{j=1}^N p_{ij}^{[M+q]}y^{[M+q]}_{j}  + { g}^{[q]}_{i} = 0, \\
    & \quad \, \, \, \forall q\in\mathcal{E}_i
\end{split}
\end{equation}
where $\mathcal{I}_i$ and $\mathcal{E}_i$ denote the set of inequality constraints and the set of equality constraints affecting agent $i$, respectively.

The Lagrangian for problem \eqref{transformed_P_i} is 
\begin{equation}\label{eq:Lagrangian}
\begin{split}
   & \mathcal{L}_i(x_i,\lambda_i, \mu_i)= f_i(x_i)\\
   &+  \sum_{m\in\mathcal{I}_i} \left \langle {\mu}_i^{[m]}, {A}_i^{[m]} x_i + y^{[m]}_{i} - \sum_{j=1}^N p_{ij}^{[m]}y^{[m]}_{j}   + { b}_i^{[m]} \right \rangle   \\
    &  + \sum_{q\in\mathcal{E}_i} \left\langle {\lambda}_i^{[q]}, E_i^{[q]} { x}_i   + y^{[M+q]}_{i} - \sum_{j=1}^N p_{ij}^{[M+q]}y^{[M+q]}_{j}  + { g}^{[q]}_{i} \right\rangle
    \end{split}
\end{equation}
where ${ \mu}_i=\text{col}(\{\mu_{i}^{[m]} \}_{m \in \mathcal{I}_i})$ and $\lambda_i=\text{col}(\{\lambda_{i}^{[q]} \}_{q \in \mathcal{E}_i})$ are the multipliers of problem \eqref{transformed_P_i} corresponding to the inequality and equality constraints, respectively. 
We note that \eqref{transformed_P_i} is a linearly constrained optimization problem. 
Thus, for fixed $y_i$'s, if $x_i(y_i)$ is an optimum of \eqref{transformed_P_i},
then there exist multipliers $\mu_i(y_i)$ and $\lambda_i(y_i)$ such that the following KKT condition holds \cite{boyd2004convex}
\begin{equation}\label{eq:KKT_system}
\begin{split}
    &\nabla_{x_i}\mathcal{L}_i\left(x_i(y_i),  \lambda_i(y_i),\mu_i(y_i)\right) = 0 \\
    & \mu_i(y_i) \geq 0\\
    & \left\langle {\mu}_i^{[m]}(y_i), {A}_i^{[m]} x_i + y^{[m]}_{i} - \sum_{j=1}^N p_{ij}^{[m]}y^{[m]}_{j}   + { b}_i^{[m]} \right\rangle = 0, \\
    &\forall m\in\mathcal{I}_i
\end{split}
\end{equation}
where
\begin{equation*}
    { y}_i=\text{col}\left(\{y^{[l]}_j\}_{j\in\mathcal{N}_i^{[l]},\forall l\in \mathcal{I}_i \cup (\{M\}\oplus \mathcal{E}_i)}
    \right)
\end{equation*}
is the collection of auxiliary variables that influence the local problem for agent $i$ in \eqref{transformed_P_i}.

We note that when solving \eqref{transformed_P_i} based on some common solvers, e.g., quadprog in MATLAB or CVX, the optimal primal and dual solutions can be simultaneously obtained.

\begin{assumption}\label{assump:licq}
For each $i\in\mathcal{V}$, there holds
\begin{itemize}
    \item[i)]  $f_i$ is $\nu$-strongly convex, i.e., there is some $\nu \geq 0$ such that
    \begin{equation*}
        f_i(y) -f_i(x) \geq \langle \nabla f_i(x), y-x \rangle + \frac{\nu}{2}\lVert y-x \rVert^2,
    \end{equation*}
    and has $\alpha$-Lipschitz gradients, i.e., there exists some $\alpha>0$ such that
        \begin{equation*}
       \lVert  \nabla f_i(y) - \nabla f_i(x) \rVert \leq {\alpha}\lVert y-x \rVert.
    \end{equation*}
    \item[ii)] the matrix whose rows are $A_i^{[m]},m\in\mathcal{I}_i$ and $E_i^{[q]}, q\in\mathcal{E}_i$
    is full row rank.
\end{itemize}
\end{assumption}

Both Assumption \ref{assump:network} and Assumption \ref{assump:licq} ii) hold under the following conditions: all $\text{mat}(A_i^T,E_i^T)$'s have full column rank, and $\mathcal{G}$ is a connected network. These conditions naturally hold in, e.g., economic dispatch \cite{yang2013consensus}.

\begin{lemma}\label{prop:duallity properties}
    Suppose Assumption \ref{assump:licq} holds. Then, for the local constrained optimization problem  in \eqref{transformed_P_i} with any given slack variables $y_i$,
\begin{itemize}
    \item[i)] feasiblility holds,
    \item[ii)] strong duality holds, and the KKT multipliers are unique.
\end{itemize}
\end{lemma}

\begin{pf}
Feasibility directly follows from Assumption \ref{assump:licq}-ii). 
With feasibility, the refined Slater's condition automatically holds if the constraints are all linear \cite[Section 5.2.3]{boyd2004convex}, which together with Assumption \ref{assump:licq}-i) yields strong duality.
In addition, the vectors
$
        \left\{ A_i^{[m]}, m\in\mathcal{I}_i; E_i^{[q]}, q\in\mathcal{E}_i \right\}
$
    are linearly independent, known as the Linear Independence Constraint Qualification (LICQ), under which the KKT multipliers for the constrained optimization problem in \eqref{transformed_P_i} are guaranteed to be unique \cite{kyparisis1985uniqueness}.  \QEDA
\end{pf}

For ease of reference, the main notations in the reformulated problem are summarized in Table \ref{table:notation}.
{\color{red}
\begin{table}[h!]
\centering
\caption{A list of main symbols in the problem reformulation.}
\begin{tabular}{ |c|c| }
\hline
$M$ & Number of inequality constraints  \\
\hline
$Q$ & Number of equality constraints  \\
\hline
\multirow{2}{*}{$\mathcal{V}^{[m]}$} & Cluster of agents affected by   \\
& the $m$-th inequality constraint  \\
\hline
\multirow{2}{*}{$\mathcal{V}^{[M+q]}$} & Cluster of agents affected by\\
& the $q$-th equality constraint  \\
\hline
$\mathcal{I}_i$ & Set of inequality constraints affecting agent $i$  \\
\hline
$\mathcal{E}_i$ & Set of equality constraints affecting agent $i$  \\
\hline
\multirow{2}{*}{$y_i^{[l]}$} & Scalar variable associated with  \\
& the $l$-th constraint at agent $i$   \\
\hline
$y^{[l]}$ & $y^{[l]}=\text{col}(\{ y^{[l]}_i \}_{i \in \mathcal{V}^{[l]}}) $ \\
\hline
$y_i$ & 
  $ { y}_i=\text{col}\left(\{y^{[l]}_j\}_{j\in\mathcal{N}_i^{[l]},\forall l\in \mathcal{I}_i \cup (\{M\}\oplus \mathcal{E}_i)} \right)$
\\
\hline
$y$ & ${ y}=\text{col}(\{y^{[l]}\}_{l=1,\dots,M+Q})$   \\
\hline
\end{tabular}
 \label{table:notation}
\end{table}}



    
    



\subsection{Perturbed function and minimization algorithm}


In this subsection, we view the problem in \eqref{transformed_P_i} as a perturbed version of the original problem \eqref{eq:centralized_problem}, defined by
\begin{equation*}
\begin{split}
    \phi_i({ y}_i) &:= \min_{ x_i} f_i(x_i) \\
      \text{s.t.} \quad   & A_{i}^{[m]} { x}_i   +  y^{[m]}_{i} - \sum_{j=1}^N p_{ij}^{[m]}y^{[m]}_{j}  + { b}^{[m]}_{i} \leq 0, \forall m\in\mathcal{I}_i \\
 &E_{i}^{[q]} { x}_i   +  y^{[M+q]}_{i} - \sum_{j=1}^N p_{ij}^{[M+q]}y^{[M+q]}_{j}  + { g}^{[q]}_{i} = 0, \\
 &\forall q\in\mathcal{E}_i.
\end{split}
\end{equation*}

Since \eqref{transformed_P_i} is always feasible under Assumption \ref{assump:licq}, $\phi_i(y_i)$ is well-defined. 
{Using the definition of $\phi(y_i)$, the optimization problem in \eqref{transformed_P}
can be equivalently expressed as 
\begin{equation}\label{eq:outer_problem}
   \min_{y} \left\{ \phi({y}) = \sum_{i=1}^N \phi_i({y}_i) \right\}.
\end{equation}}

Upon substituting the optimal perturbation $y^*$ into  \eqref{transformed_P}, the original problem \eqref{eq:centralized_problem} can be solved. 

{
\begin{remark}
    The rank condition in Assumption \ref{assump:licq} can be relaxed. Indeed, in the case where Assumption \ref{assump:licq}-ii) fails, one can consider a relaxed version of the perturbed problem:
    \begin{equation*}
\begin{split}
    \phi'_i({ y}_i) &:= \min_{ x_i, \rho_i\geq 0} f_i(x_i)+\omega\rho_i, \\
         \text{s.t.} \quad &A_{i}^{[m]} { x}_i   +  y^{[m]}_{i} - \sum_{j=1}^N p_{ij}^{[m]}y^{[m]}_{j}  + { b}^{[m]}_{i} \leq \rho_i,  \forall m\in\mathcal{I}_i \\
 &  E_{i}^{[q]} { x}_i   +  y^{[M+q]}_{i} - \sum_{j=1}^N p_{ij}^{[M+q]}y^{[M+q]}_{j}  + { g}^{[q]}_{i} =0, \\
 &\forall q\in\mathcal{E}_i
\end{split}
\end{equation*}
where $\omega$ is a positive scalar. If \eqref{eq:centralized_problem} is feasible and $\omega$ is sufficiently large, then 
there must exist some $y_i$ such that $x_i^*$ and $\rho_i=0$, $i=1,\cdots, N$ are optimal to the perturbed problem $\min_y \sum_{i=1}^n \phi'_i({ y}_i)$. Thus, by solving $\min_y \sum_{i=1}^n  \phi'_i({ y}_i)$, one attains an optimal solution to \eqref{eq:centralized_problem} \cite[Proposition III.3]{notarnicola2019constraint}.
Therefore, the algorithms presented next also apply to this setup, with minor modifications on the subproblem.
However, for this relaxed problem, the projection of its feasible set onto the space of $x_i$ is not identical to that of the original problem,
i.e., the iterates during implementation are no longer  guaranteed to be feasible to \eqref{eq:centralized_problem}.
\end{remark} }

It is noteworthy that there is no consensus constraint when minimizing $\phi(y)$. As a result, further decomposition of this problem is not required, which is a key difference between the proposed method with existing literature \cite{notarnicola2019constraint}. Next, we provide explicit expressions for the gradients of $\phi(y)$ in the following lemma.

\begin{lemma}\label{prop:global_subgradient}
    Suppose Assumptions \ref{assump:network} and \ref{assump:licq} hold with $\nu>0$. There holds
    \begin{itemize}
        \item[i)] $\phi({y})$ is convex and differentiable,
        \item[ii)] the gradients of $\phi({y})$ can be computed as 
        \begin{equation}\label{eq:gradient}
        \begin{split}
            \nabla_{y^{[m]}} \phi(y) & =(I-P^{[m]})^T \text{col}( \{\mu_i^{[m]}(y_i)\}_{i\in\mathcal{V}^{[m]}} )  \\
 \nabla_{y^{[M+q]}} \phi(y)      &= (I-P^{[M+q]})^T \text{col}( \{\lambda_i^{[q]}(y_i)\}_{i\in\mathcal{V}^{[M+q]}} )      
                    \end{split}
        \end{equation}
        where $\mu_i^{[m]}(y_i)$ and $\lambda_i^{[q]}(y_i)$ denote the KKT multipliers associated with the corresponding inequality and equality constraints in \eqref{transformed_P_i},
        \item[iii)] $\phi({y})$ has $\alpha_\phi$-Lipschitz continuous gradients for some finite $\alpha_\phi$. 
    \end{itemize}
\end{lemma}

\begin{pf}
i) The convexity of $\phi$ can be proved by following \cite[Section 5.6.1]{boyd2004convex}. For completeness, a proof is included. Let $x$ and $x'$ be the optimal primal solutions of \eqref{eq:outer_problem} under $y$ and $y'$. Given any $\lambda\in[0,1]$, $\phi(\lambda y+(1-\lambda)y')$ is well-defined, and there holds
\begin{equation*}
\begin{split}
    &\lambda \phi(y)+(1-\lambda)\phi(y')\\
    &= \lambda\sum_{i=1}^Nf_i( x )+ (1-\lambda)\sum_{i=1}^Nf_i( x' ) \\
    &\geq  \sum_{i=1}^Nf_i( \lambda x+(1-\lambda)x' )
\end{split}
\end{equation*}
where the inequality is due to the convexity of $f_i$. In addition, $\lambda x+ (1-\lambda) x'$ is feasible to the optimization problem corresponding to $\phi(\lambda y+ (1-\lambda) y')$ due to linearity of the constraints. Thus, by optimality, we obtain
\begin{equation*}
\begin{split}
    \lambda \phi(y)+(1-\lambda)\phi(y')
    \geq \phi(\lambda y+(1-\lambda)y').
\end{split}
\end{equation*}
According to Lemma \ref{prop:duallity properties}, the KKT multipliers are unique for any auxiliary variable $y$, and therefore $\phi(y)$ is globally differentiable \cite[Corollary 7.3.1]{florenzano2001finite}.

   ii) Because strong duality holds, we obtain from sensitivity analysis \cite[Section 5.6.2]{boyd2004convex} that
\begin{equation*}
\begin{split}
    &\phi_i({ y}_i') \\
    &\geq \phi_i({ y}_i) + \sum_{m\in\mathcal{I}_i}\left\langle {\mu}_i^{[m]}(y_i), 
     (I-P^{[m]})_{i}\left(({ y}')^{[m]} - y^{[m]}\right) \right\rangle  \\
    & \quad +  \sum_{q\in\mathcal{E}_i} \left\langle {\lambda}_i^{[q]}(y_i),  (I-P^{[M+q]})_{i}\left(({ y}')^{[M+q]} - y^{[M+q]}\right) \right\rangle \\
    & =  \phi_i({ y}_i) + \sum_{m\in\mathcal{I}_i} \left\langle  (I-P^{[m]})_{i}^T{\mu}_i^{[m]}(y_i), 
   ({ y}')^{[m]}- y^{[m]} \right\rangle  \\
   & \quad +  \sum_{q\in\mathcal{E}_i} \left\langle (I-P^{[M+q]})_{i}^T {\lambda}_i^{[q]}(y_i),  ({ y}')^{[M+q]}-y^{[M+q]} \right\rangle
    \end{split}
\end{equation*}
where $(\cdot)_i$ represents the $i$-th row of a matrix.
This indicates that 
\begin{equation*}
  \nabla_{y^{[m]}} \phi_i({ y}_i)
  =  
  (I-P^{[m]})_i^T
  \mu_i^{[m]}(y_i), \forall  m\in\mathcal{I}_i
\end{equation*}
and
\begin{equation*}
  \nabla_{y^{[M+q]}} \phi_i({ y}_i) 
  =  
  (I-P^{[M+q]})_i^T
  \lambda_i^{[q]}(y_i), \forall q\in\mathcal{E}_i
\end{equation*}
where $\nabla_{y^{[l]}}\phi_i$ 
denotes the block coordinate gradient with respect to $y^{[l]}$.
{Note that $\phi({y}) = \sum_{i=1}^N \phi_i({y}_i)$ and each $y_i$ is partially coupled with each other.
For each coordinate block indexed by $m$, $\nabla_{y^{[m]}} \phi({ y})$ is contributed only by $ \mu_i^{[m]}(y_i)(I-P^{[m]})_i^T, i\in\mathcal{V}^{[m]}$, that is, $\forall m=1,\dots, M$
\begin{equation*}
\begin{split}
     &\nabla_{y^{[m]}} \phi({ y})
     =  \sum_{i\in\mathcal{V}^{[m]}} \mu_i^{[m]}(y_i)(I-P^{[m]})_i^T \\
       &
       = (I-P^{[m]})^T \text{col}( \{\mu_i^{[m]}(y_i)\}_{i\in\mathcal{V}^{[m]}} )
\end{split}
\end{equation*}}
and, similarly, $\forall q = 1,\dots, Q$
\begin{equation}\label{eq:subgradient_z}
    \nabla_{y^{[M+q]}} \phi({ y}) =
    (I-P^{[M+q]})^T \text{col}( \{\lambda_i^{[q]}(y_i)\}_{i\in\mathcal{V}^{[M+q]}} ).
\end{equation}

iii) Because each $f_i$ is strongly convex, the second-order sufficient conditions (SOSC) automatically hold, which implies that the KKT multipliers are Lipschitz continuous with finite positive parameters \cite{subotic2021quantitative}. Based on the formulas of gradients given in ii), we obtain that $\phi$ has Lipschitz continuous gradients with a finite positive parameter. \QEDA
\end{pf}

{
We have converted the original problem represented in \eqref{eq:centralized_problem} into the form \eqref{eq:outer_problem}. Consequently, finding a solution for \eqref{eq:centralized_problem} is equivalent to solving \eqref{eq:outer_problem}. To accomplish this, it becomes essential to compute the block coordinate gradient of $\nabla_{y^{[l]}} \phi({ y})$ in a distributed manner.
Lemma \ref{prop:global_subgradient} indicates that each agent $i\in\mathcal{V}^{[l]}$ can collaboratively compute $\nabla_{y^{[l]}} \phi({ y})$, $l=1,\dots, M+Q$ by communicating with their immediate neighbors about the KKT multipliers. This is partially due to the fact that the closed-form expression of the block coordinate gradient involves a linear mapping $(I-P^{[m]})^T$, which has been intentionally designed to be compatible with the communication network.}
Based on this fact, the accelerated dual averaging method \cite{cohen2018acceleration} can be implemented to solve \eqref{eq:outer_problem} in a distributed manner.

The overall algorithm is summarized in Algorithm \ref{algo:DO}. 
{Each agent updates two sequences of auxiliary variables: $y_i$ and $\hat{y}_i$, according to the following rules.
In step 5, each agent communicates with its neighbors to collect $y_j$ variables to formulate a local problem in the form of \eqref{transformed_P_i}. After locally identifying the Lagrangian multipliers for \eqref{transformed_P_i}, each agent exchanges the multipliers with its neighbors in order to compute the block coordinate gradient of $\phi(y)$, as outlined in steps 6-8. Finally, each agent follows the accelerated dual averaging to update $y_i$ in step 12. We remark that the parameters $\gamma_t$ and $\Gamma_t$ should be properly chosen to accelerate the convergence \cite{cohen2018acceleration}, as outlined in Theorem \ref{thm:convergence_GD}.
} 


\begin{algorithm}[tb]
	\caption{Violation-free distributed accelerated dual averaging}
	\label{algo:DO}
	\begin{algorithmic}[1]
		\STATE {\bfseries input:}  arbitrary variable $y$, parameters $\{\gamma_t\}_{t\geq 0}$ and $\{ \Gamma_t = \sum_{\tau=1}^t \gamma_\tau \}_{t\geq 0}$
  	\STATE {\bfseries output}:
        $x_i(\hat{y}_i), i=1,\dots,n$
  \FOR{$t=1,2,\dots $}
      \STATE for each agent $i\in \{1,\dots, N\}$:
			\STATE collect $y_j^{[l]}$ from each $ j\in\mathcal{N}_{i}^{[l]}$,  $\forall l\in \mathcal{I}_i \cup (\{M\}\oplus \mathcal{E}_i)$
		\STATE compute $x_i(y_i)$ and the multipliers $\mu_i(y_i)$ and $\lambda_i(y_i)$  by solving \eqref{transformed_P_i}
  \STATE collect $\mu_j^{[m]}(y_j)$ from each $j\in\mathcal{N}_{i}^{[m]}$, $ \forall m\in \mathcal{I}_i$ and $\lambda_i^{[q]}(y_i)$ from each $j\in\mathcal{N}_{i}^{[M+q]}$,$\forall q\in\mathcal{E}_i$
  \STATE compute $\nabla_{y_i^{[l]}} \phi(y), \forall l\in \mathcal{I}_i \cup (\{M\}\oplus \mathcal{E}_i)$ 
  according to \eqref{eq:gradient}

  \IF{$t=1$}
  \STATE set $\hat{y}_i^{[l]}= z_i^{[l]}= -\gamma_1\nabla_{y_i^{[l]}}\phi(y), \forall l\in \mathcal{I}_i \cup (\{M\}\oplus \mathcal{E}_i)$
  \ELSE
  \STATE update $y_i^{[l]}, \forall l\in \mathcal{I}_i \cup (\{M\}\oplus \mathcal{E}_i)$ by
    \begin{equation*}
      \begin{split}
          y_i^{[l]} &\leftarrow \left(1-\frac{\gamma_t}{\Gamma_t}\right) \hat{y}_i^{[l]} + \frac{\gamma_t}{\Gamma_t}  z_i^{[l]} \\
          z_i^{[l]}  &\leftarrow  z_i^{[l]} - \gamma_t \nabla_{y_i^{[l]}} \phi(y) \\
          \hat{y}_i^{[l]} &\leftarrow \left(1-\frac{\gamma_t}{\Gamma_t}\right) \hat{y}_i^{[l]} + \frac{\gamma_t}{\Gamma_t}z_i^{[l]}
      \end{split}
  \end{equation*}
    \ENDIF
		\ENDFOR
	\end{algorithmic}
\end{algorithm}

\section{Convergence analysis}
\label{sec:4}
\subsection{Computing Lipschitz constants}

Before establishing conditions under which Algorithm \ref{algo:DO} converges, we quantify the Lipschitz constant $\alpha_\phi$ of the gradients of $\phi(y)$.
In particular, we present the Lipschitz constants of $\mu_i(y_i)$ and $\lambda_i(y_i)$, based on which the Lipschitz constant of $\nabla\phi(y)$ can be computed.

Under Assumption \ref{assump:licq} with $\nu>0$, \eqref{transformed_P} has a unique regular optimizer for any $y^{[l]}, l=1,\dots, M+Q$. In addition, the KKT multipliers are globally continuous with respect to $y$ \cite[Theorem 2]{subotic2021quantitative}. 

\begin{proposition}
    Suppose Assumptions \ref{assump:network} and \ref{assump:licq} hold with $\nu>0$. The Lipschitz constant of the gradients of $\phi(y)$ can be computed as
\begin{equation}\label{eq:smoothness_param}
\begin{split}
    &\alpha_{\phi} \leq \\
    & (\max_{i} \alpha_{(\lambda_i,\mu_i)} ) \left(\max_{l\in\{1,\dots,M+Q\} }\lVert I-P^{[l]} \rVert  
 \sqrt{\lvert \mathcal{V}^{[l]} \rvert} \right)  \sqrt{M+Q}
\end{split}
\end{equation} 
where $\alpha_{(\lambda_i,\mu_i)}$ denotes the Lipschitz constant of $\mu_i(y_i)$ and $\lambda_i(y_i)$, given by
\begin{equation*}
\begin{split}
    \alpha_{(\lambda_i,\mu_i)}\leq 
   \left(\max_{l \in\mathcal{I}_i \cup (\{M\}\oplus \mathcal{E}_i)}  \lVert I-P^{[l]} \rVert  \right) \sqrt{\frac{{\alpha}}{\lambda_{\text{min}}(\hat{\mathcal{B}}_i^{T} \hat{\mathcal{B}}_i)}} 
\end{split}
\end{equation*}
where $\alpha$ denotes the Lipschitz constant of $\nabla f_i$,  $\hat{\mathcal{B}}_i=\text{mat}(\hat{E}_i^{T}, 
        \hat{A}_i^T)$, and $\text{mat}(A,B)=\begin{bmatrix}
    A, B
\end{bmatrix}$.
\end{proposition}

\begin{pf}
Recall the Lagrangian in \eqref{eq:Lagrangian} and consider the KKT system
\begin{equation*}
    \begin{split}
&\mathcal{F}_i(x_i,\lambda_i,\mu_i,y_i):=\\
&\begin{bmatrix}
        \nabla_{x_i} \mathcal{L}_i(x_i, \lambda_i, \mu_i) \\
        \hline
        \hat{E}_i x_i + \hat{g}_i + \text{blkdiag}(\{ (I-P^{[M+q]})_i\}_{q\in\mathcal{E}_i}) y_i^{\mathcal{E}_i} \\
        \hline
        \text{diag}(\{\mu_i^{[m]}\}_{m\in\mathcal{I}_i}) \times \\
        \left(    \hat{A}_i x_i + \hat{b}_i + \text{blkdiag}(\{ (I-P^{[m]})_i\}_{m\in\mathcal{I}_i}) y_i^{\mathcal{I}_i} \right)
    \end{bmatrix} = 0 
    \end{split}
\end{equation*}
where
\begin{equation*}
{y}_i^{\mathcal{I}_i}=\text{col}(\{y^{[m]}\}_{m\in\mathcal{I}_i}) \,\, \text{and} \,\, { y}_i^{\mathcal{E}_i}=\text{col}(\{y^{[M+q]}\}_{q\in\mathcal{E}_i}),
\end{equation*}
and 
\begin{equation*}
\hat{E}_i = [ E_i^{[1]}; \dots; E_i^{\lvert \mathcal{E}_i \rvert } ] \,\, \text{and} \,\,  \hat{g}_i = \text{col} (\{ g_i^{[q]} \}_{q\in\mathcal{E}_i}),
\end{equation*}
\begin{equation*}
  \hat{A}_i = [ A_i^{[1]}; \dots; A_i^{\lvert  \mathcal{I}_i \rvert }  ]  \,\, \text{and} \,\, \hat{b}_i = \text{col}( \{ b_i^{[m]} \}_{m\in\mathcal{I}_i}).
\end{equation*}

Next, to ease notation, we omit the argument from any map that solely depends on $y$. 

We bound the Lipschitz constants for functions $\mu_i^{[m]}$ and $\lambda_i^{[q]}$ by following the approach in \cite[Proposition 4]{subotic2021quantitative}. Denote by \begin{equation*}
  \mathbf{I}_i:=\left\{ m:  A_i^{[m]}x_i + y_i^{[m]}-\sum_{j \in \mathcal{N}_{i}^{[m]}} p_{ij}^{[m]} y_j^{[m]} + b_i^{[m]} = 0  \right\}
\end{equation*} 
 the set of active inequality constraints for agent $i$, and by $\overline{\mathbf{I}}_i$ the set of inactive constraints for $i$. Without loss of generality, we assume the first $\lvert  \mathbf{I}_i \rvert$ inequality constraints are active. For simplicity, the inequality constraints are required to be strictly active, that is, the KKT multiplier associated with the constraint should be strictly positive. As a consequence, the KKT multipliers are differentiable with respect to $y_i$ \cite{subotic2021quantitative}. However, such bound can be generalized to the case with weakly active inequality constraints by following \cite{jittorntrum1978sequential}.


By the implicit function theorem \cite[Appendix A]{bertsekas2016nonlinear} and by organizing the matrices according to $
    (x_i,[\lambda_i, \mu_i^{{\bf I}_i} ], \mu_i^{\overline{\mathbf{I}}_i})$, 
it holds that
\begin{equation}\label{eq:implicit_function_thm}
    \nabla_{y_i} (x_i,\lambda_i,\mu_i) =   -(\nabla_{x_i,\lambda_i,\mu_i} \mathcal{F}_i )^{-1} \nabla_{y_i} \mathcal{F}_i
\end{equation}
where
\begin{equation*}
    \nabla_{x_i,\lambda_i, \mu_i}\mathcal{F}_i:=\begin{bmatrix}
        \nabla_{x_ix_i}^2\mathcal{L}_i & \mathcal{B}_i & \hat{A}_i^{\overline{\bf I}_i} \\
        \mathcal{D}_{\mathbf{I}_i}\mathcal{B}_i^T & 0 & 0\\
        0 & 0 & \text{diag}( \mathcal{H}_i^{\overline{\mathbf{I}}_i}  )
    \end{bmatrix},
\end{equation*}
\begin{equation*}
    \nabla_{y_i} \mathcal{F}_i := \begin{bmatrix}
        0 \\ \mathcal{D}_{\mathbf{I}_i}\mathcal{G}_i \\ 0
    \end{bmatrix}, \,\, 
    \mathcal{B}_i =
       \text{mat} (\hat{E}_i^{T}, 
        \hat{A}_i^{\mathbf{I}_i}),
\end{equation*}
\begin{equation*}
  \mathcal{D}_{{\bf I}_i} = \text{diag}([{\bf 1}_{\lvert\mathcal{E}_i\rvert}; \text{col}(\{ \mu_i^{[m]}\}_{m\in{\mathbf{I}}_i})]),
\end{equation*}
\begin{equation*}
 \hat{A}_i^{{\bf I}_i} = 
         \text{mat}(\{ (A_i^{[m]})^T \}_{m\in  {\bf I}_i}), \,\,   \hat{A}_i^{\overline{\bf I}_i} = 
        \text{mat} (\{ (A_i^{[m]})^T \}_{m\in  \overline{\bf I}_i}) ,
\end{equation*}
\begin{equation*}
    \mathcal{G}_i = \text{blkdiag}(\{ (I-P^{[M+q]})_i\}_{q\in\mathcal{E}_i},\{ (I-P^{[m]})_i\}_{m\in\mathbf{I}_i}),
\end{equation*}
\begin{equation*}
    \mathcal{H}_i =   
          \hat{A}_i x_i + \hat{b}_i + \text{blkdiag}(\{ (I-P^{[m]})_i\}_{m\in\mathcal{I}_i}) y_i^{\mathcal{I}_i},
\end{equation*}
$\mathcal{H}_i^{\overline{\mathbf{I}}_i} = \text{col}(\{ \mathcal{H}_i^{[m]} \}_{m\in  \overline{\bf I}_i})$, and $\nabla_{x_ix_i}^2$ represents the Hessian operator.
It can be verified that
\begin{equation*}
    (\nabla_{x_i,\lambda_i, \mu_i}\mathcal{F}_i)^{-1}=\begin{bmatrix}
       \mathcal{M}_i^{-1} & *  \\
        0 & * 
    \end{bmatrix}\,\,  \text{with} \,\, \mathcal{M}_i =  \begin{bmatrix}
        \nabla_{x_ix_i}^2\mathcal{L}_i & \mathcal{B}_i \\
           \mathcal{D}_{\mathbf{I}_i}\mathcal{B}_i^T & 0 
    \end{bmatrix}
\end{equation*}
where the symbol $*$ represents non-zero components that are irrelevant to the context.
Therefore, for the inactive inequality constraints, it holds that 
$
    \nabla_{y_i} \mu_{\overline{\bf I}_i} = 0.
$
And for equality and active inequality constraints, we obtain from the formula of the inverse of $2\times 2$ block matrices \cite{lu2002inverses} that
\begin{equation*}
    \nabla_{y_i}\begin{bmatrix}
        \lambda_i \\
        \mu_i^{\mathbf{I}_i} 
    \end{bmatrix} = (\mathcal{M}_i^{-1})_{22}\mathcal{D}_{\mathbf{I}_i}\mathcal{G}_i = (\mathcal{B}_i^T (\nabla_{xx}^2 \mathcal{L}_i)^{-1}\mathcal{B}_i )^{-1}
    \mathcal{G}_i
\end{equation*}
where $(\mathcal{M}_i^{-1})_{22}$ denotes the $(2,2)$-th entry of $\mathcal{M}_i^{-1}$. This indicates that the Lipschitz constants of $\mu_i$ and $\lambda_i$
\begin{equation*}
\begin{split}
   & \alpha_{(\lambda_i,\mu_i)} \leq \lVert 
    (\mathcal{B}_i^T (\nabla_{xx}^2 \mathcal{L}_i)^{-1}\mathcal{B}_i )^{-1}
    \mathcal{G}_i\rVert  \\
    & \leq  \lVert 
    (\mathcal{B}_i^T (\nabla_{xx}^2 \mathcal{L}_i)^{-1}\mathcal{B}_i )^{-1}
   \rVert \lVert 
    \mathcal{G}_i\rVert \\
    & \leq \lVert 
    \mathcal{G}_i\rVert  \sqrt{\frac{{\alpha}}{\lambda_{\text{min}}(\hat{\mathcal{B}}_i^{T} \hat{\mathcal{B}}_i)}} \\
    & \leq  \left(\max_{l \in\mathcal{I}_i \cup (\{M\}\oplus \mathcal{E}_i)}  \lVert I-P^{[l]} \rVert  \right)  \sqrt{\frac{{\alpha}}{\lambda_{\text{min}}(\hat{\mathcal{B}}_i^{T} \hat{\mathcal{B}}_i)}} 
\end{split}
\end{equation*}
where the third and fourth inequality follows from
\begin{equation*}
    \lVert 
    \mathcal{B}_i^T (\nabla_{xx}^2 \mathcal{L}_i)^{-1}\mathcal{B}_i 
   \rVert \geq \sqrt{\frac{\lambda_{\text{min}}({\mathcal{B}}_i^{T} {\mathcal{B}}_i)}{{\alpha}}} \geq  \sqrt{\frac{\lambda_{\text{min}}(\hat{\mathcal{B}}_i^{T} \hat{\mathcal{B}}_i)}{{\alpha}}}
\end{equation*}
and
\begin{equation*}
    \lVert \text{blkdiag}(A,B) \rVert = \text{max}\{ \lVert A \rVert, \lVert B \rVert \}\,\, \text{and} \,\, \lVert (A)_i \rVert \leq  \lVert A \rVert,
\end{equation*}
respectively.
We note that the above bound of $\alpha_{(\lambda_i,\mu_i)}$ holds for uncertain $\mathbf{I}_i$.

Given $y$ and $y'$, we have
\begin{equation*}
\begin{split}
    & \lVert \text{col}( \{\mu_i^{[m]}(y_i)\}_{i\in\mathcal{V}^{[m]}} ) - \text{col}( \{\mu_i^{[m]}(y_i')\}_{i\in\mathcal{V}^{[m]}} )    \rVert  \\ 
    & = \sqrt{ \sum_{i\in\mathcal{V}^{[m]}}\lVert  \mu_i^{[m]}(y_i)  -  \mu_i^{[m]}(y_i')    \rVert^2  }  \\
    & \leq (\max_{i} \alpha_{(\lambda_i,\mu_i)} ) \sqrt{ \sum_{i\in\mathcal{V}^{[m]}}\lVert  y_i  -  y_i'    \rVert^2  } \\
    &= (\max_{i} \alpha_{(\lambda_i,\mu_i)} )\sqrt{\lvert \mathcal{V}^{[m]} \rvert}  \lVert  y  -  y'    \rVert.
\end{split}
\end{equation*}
This in conjunction with Lemma \ref{prop:global_subgradient} leads to
\begin{equation*}
\begin{split}
    & \lVert  \nabla_{y^{[m]}}\phi(y) - \nabla_{y^{[m]}}\phi(y') \rVert  \\
    & \leq  \lVert I-P^{[m]} \rVert  
(\max_{i} \alpha_{(\lambda_i,\mu_i)} ) \sqrt{\lvert \mathcal{V}^{[m]} \rvert}  \lVert  y  -  y'    \rVert.     
\end{split}
\end{equation*}
Furthermore, 
\begin{equation*}
\begin{split}
    & \lVert  \nabla_{y}\phi(y) - \nabla_{y}\phi(y') \rVert  \\
    & = \sqrt{ \sum_{l=1}^{M+Q}\lVert  \nabla_{y^{[l]}}\phi(y) - \nabla_{y^{[l]}}\phi(y') \rVert^2   } \\
    & \leq (\max_{i} \alpha_{(\lambda_i,\mu_i)} ) \left(\max_{l\in\{1,\dots,M+Q\} }\lVert I-P^{[l]} \rVert  
 \sqrt{\lvert \mathcal{V}^{[l]} \rvert} \right) \\
 & \quad \times \sqrt{M+Q}\lVert  y  -  y'    \rVert.     
\end{split}
\end{equation*}
The proof is complete.
\QEDA
\end{pf}

\subsection{Rate of convergence}

\begin{thm}\label{thm:convergence_GD}
    Suppose the premise given in Lemma \ref{prop:global_subgradient} holds. If the parameter $\gamma_t$ in Algorithm \ref{algo:DO} is set as $\gamma(t+1)$ where
    $
        \gamma \leq  1/{(2\alpha_{\phi})},
    $
    then, for all $t\geq 2$, the solution to the problem in \eqref{transformed_P} with auxiliary variable $ y=\hat{y}^{(t)}$ is a feasible solution to \eqref{eq:centralized_problem} and satisfies
    \begin{equation*}
        \sum_{i=1}^N f_i({x}_i(\hat{y}_i^{(t)})) - f_i(x_i^{*}) \leq  \frac{\lVert {y}^{(0)}-y^* \rVert^2}{2\Gamma_t} = \frac{2 \alpha_\phi\lVert {y}^{(0)}-y^* \rVert^2}{t(t+3)}
    \end{equation*}
    where $y^{(0)}$ denotes the initial variable and $y^*$ is any (fixed) minimizer of \eqref{eq:outer_problem}.
\end{thm}
\begin{pf}
The proof is a consequence of \cite[Theorem 3.4]{cohen2018acceleration} and omitted here for brevity. \QEDA

\end{pf}

\section{Extension to non-strongly convex case}
\label{sec:5}

In this section, we relax the strong convexity assumption on $f_i$. 

For general convex $f_i$, following the same proof of Lemma \ref{prop:duallity properties} $\phi(y)$ remains differentiable. However, the KKT multipliers as well as the gradients of $\phi(y)$ are not necessarily Lipschitz continuous. In addition, the minimization problem in \eqref{eq:outer_problem} is unconstrained and has noncompact solution set. 
To guarantee the boundedness of $\nabla\phi(y)$, we introduce a compact constraint
\begin{equation*}
    \mathcal{Y} := [-C, C]^{\sum_{l=1}^{M+Q}\lvert \mathcal{V}^{[l]} \rvert}  
\end{equation*}
to $y$, where $C$ is a positive constant. Note that  $\mathcal{Y}$ imposes constraint on each individual coordinate of $y$, which facilitates projection in local updates of $y_i^{[m]}$. If $C$ is sufficiently large, this constraint will not change the optimal solution of the problem in \eqref{eq:outer_problem}. Upon incorporating the constraints, we arrive at the following optimization problem:
\begin{equation}\label{eq:outer_problem_constrained}
   \min_{y\in\mathcal{Y}} \phi({y}).
\end{equation}

\begin{proposition}\label{prop:identical_cost}
    Let $y^*$ be an optimal solution of problem \eqref{eq:outer_problem} and $C$ is a positive constant such that $y^*\in\mathcal{Y}$. Then, $y^*$ is also an optimal solution of problem \eqref{eq:outer_problem_constrained} and the two problems have identical optimal cost.
\end{proposition}
\begin{pf}
    Since problem in \eqref{eq:outer_problem} is unconstrained, it has lower or identical optimal cost than that in \eqref{eq:outer_problem_constrained}. Because $y^*\in\mathcal{Y}$, they lead to a cost of \eqref{eq:outer_problem_constrained} that is equal to the optimal cost \eqref{eq:outer_problem}, implying that $y^*$ is also an optimal solution of \eqref{eq:outer_problem_constrained} and the two problems have identical optimal cost. \QEDA
\end{pf}

The compactness of $\mathcal{Y}$, together with the convexity of $\phi(y)$, leads to the following lemma, whose proof can be found in \cite[Appendix B]{bertsekas2016nonlinear}.

\begin{lemma}\label{lem:boundedness}
   Under Assumption \ref{assump:licq}, the gradients of $\phi(y)$ in \eqref{eq:outer_problem_constrained} are bounded from above, i.e., $\lVert \nabla \phi(y) \rVert \leq G$.
\end{lemma}

Based on Lemma \ref{lem:boundedness}, gradient descent with diminishing step-size can be used to solve \eqref{eq:outer_problem_constrained} with convergence guarantee \cite[Theorem 8.30]{beck2017first}. The algorithm is summarized in Algorithm \ref{algo:DO_nonSC}.

\begin{algorithm}[tb]
	\caption{Violation-free distributed gradient method}
	\label{algo:DO_nonSC}
	\begin{algorithmic}[1]
		\STATE {\bfseries input:}  arbitrary variable $y$, and parameter $\{\gamma_t\}_{t\geq 0}$ 
    	\STATE {\bfseries output}:
        $x_i({y}_i), i=1,\dots,n$
  \FOR{$t=1,2,\dots $}
      \STATE for each agent $i\in \{1,\dots, N\}$:
			\STATE collect $y_j^{[l]}$ from each $ j\in\mathcal{N}_{i}^{[l]}$,  $\forall l\in \mathcal{I}_i \cup (\{M\}\oplus \mathcal{E}_i)$
		\STATE compute $x_i(y_i)$ and the multipliers $\mu_i(y_i)$ and $\lambda_i(y_i)$  by solving \eqref{transformed_P_i}
  \STATE collect $\mu_j^{[m]}(y_j)$ from each $j\in\mathcal{N}_{i}^{[m]}$, $\forall m\in\mathcal{I}_i$, and $\lambda_i^{[q]}(y_i)$ from each $j\in\mathcal{N}_{i}^{[M+q]}$,$\forall q\in\mathcal{E}_i$
  \STATE compute $\nabla_{y_i^{[l]}} \phi(y), \forall l\in \mathcal{I}_i \cup (\{M\}\oplus \mathcal{E}_i)$ 
  according to \eqref{eq:gradient}
  \STATE update $y_i^{[l]}, \forall l\in \mathcal{I}_i \cup (\{M\}\oplus \mathcal{E}_i)$
  by 
    \begin{equation*}
      \begin{split}
          y_i^{[l]} &\leftarrow \text{Proj}_{[-C,C]} ({y}_i^{[l]} -{\gamma_t}  \nabla_{y_i^{[l]}} \phi(y)) \\
      \end{split}
  \end{equation*}
  
		\ENDFOR
	\end{algorithmic}
\end{algorithm}

{
\begin{thm}
    Suppose Assumptions \ref{assump:network} and \ref{assump:licq} hold with $\nu=0$. If the parameter $\gamma_t$ in Algorithm \ref{algo:DO_nonSC} is set as $\frac{\sqrt{2\Theta}}{G\sqrt{t+1}}$ where $\Theta=\max_{x,y\in 
    \mathcal{Y}}\lVert x-y \rVert^2/2$,
    then, for all $t\geq 2$, the solution to the problem in \eqref{transformed_P} is a feasible solution to \eqref{eq:centralized_problem} and satisfies
    \begin{equation*}
                \min_{\tau=0,1,\dots,t} \sum_{i=1}^N f_i({x}_i({y}_i^{(\tau)})) 
        - f_i(x_i^{*}) \leq \frac{2(1+\log(3))G\sqrt{2\Theta}}{\sqrt{t+2}}.
    \end{equation*}
\end{thm}}


\section{Numerical experiment}
\label{sec:6}

In this section, we apply the proposed algorithm to solve a multi-agent consensus problem subject to two coupling state constraints.

\subsection{Problem setup}

Consider a multi-agent system with $N=7$ agents. They are connected via a line graph. For each agent $i$, the state is denoted as $z_i\in\mathbb{R}^2$ and the system dynamics is 
$
    \dot{z}_i = x_i.
$
The states are initialized as $
    z_i(0) = (2\cos(2\pi i/7)+2, 2\sin(2\pi i/7)+1), i=1,\dots, N.
$
Define $z=\text{col}(\{ z_i\}_{i=1,\dots,7})$ and $x=\text{col}(\{ x_i\}_{i=1,\dots,7})$.
The agents seek consensus under two coupling state constraints given by
\begin{equation}\label{eq:CM_constraint}
\begin{split}
    \mathcal{Z} := \Big\{&z\in\mathbb{R}^{14}, 
        g_1(\text{col}(\{z_i\}_{i=1,\dots,4}))\geq 0, \\
       & g_2(\text{col}(\{z_i\}_{i=4,\dots,7}))\geq 0 \Big\}
\end{split}
\end{equation}
where
\begin{equation}\label{eq:constraint_functions}
\begin{split}
   & g_1(\text{col}(\{z_i\}_{i=1,\dots,4})) =  4-\sum_{i=1}^4z_i^Tz_i  \\
   & g_2(\text{col}(\{z_i\}_{i=4,\dots,7})) = 16-\sum_{i=4}^7 (z_i-[2;2])^T  (z_i-[2;2]) .
\end{split}
\end{equation}
To solve this constrained control problem, we use the consensus protocol $ 
    x_{nom,i} = \sum_{j\in\mathcal{N}_i}(z_j -z_i)$ and follow the CBF approach \cite{xu2015robustness,ames2016control} to handle the constraints. To proceed, we present the definition of CBF particularly for the dynamical system 
    \begin{equation}\label{eq:system_dynamics}
      \dot{z}={x}.
    \end{equation}
\begin{definition}[CBF] 
Consider the safety set 
$
    \mathcal{Z}$ and the dynamical system given by \eqref{eq:CM_constraint} and \eqref{eq:system_dynamics}, respectively. The functions in \eqref{eq:constraint_functions} are CBFs if there exists a locally Lipschitz, strictly increasing function $\alpha(\cdot)$ with $\alpha(0)=0$ such that $\forall z$, $\exists { x}$
    \begin{equation}\label{eq:const1}
    \begin{split}
   &\left\langle \nabla g_{1}(\text{col}(\{z_i\}_{i=1,\dots,4})), \text{col}(\{x_i\}_{i=1,\dots,4})  \right\rangle  \\
   &+ \alpha (g_1(\text{col}(\{z_i\}_{i=1,\dots,4}))) \geq 0, 
       \end{split}
\end{equation}
and
\begin{equation}\label{eq:const2}
   \begin{split}
    &\left\langle \nabla g_{2}(\text{col}(\{z_i\}_{i=4,\dots,7})), \text{col}(\{x_i\}_{i=4,\dots,7})  \right\rangle   \\
    &+ \alpha (g_2(\text{col}(\{z_i\}_{i=4,\dots,7}))) \geq 0.
        \end{split}
\end{equation}
\end{definition}
It has been shown that any locally Lipschitz ${x}$ that fulfils the CBF constraint in \eqref{eq:const1} and \eqref{eq:const2} makes the set $\mathcal{Z}$ forward invariant, and if $\mathcal{Z}$ is compact, asymptotically stable \cite{xu2015robustness}.
Following this fact and setting $\alpha$ as the identity function, for any given states $z$, an optimization-based controller can be constructed as
\begin{equation}\label{eq:CBF_constraint}
\begin{split}
    &\min_{x_1,\dots, x_N} \sum_{i=1}^N \frac{1}{2}\lVert x_i -x_{nom,i} \rVert^2 \\
     \text{s.t.} \quad 
    & \text{constraints in \eqref{eq:const1} and \eqref{eq:const2}}. 
\end{split}
\end{equation}

\subsection{Performance on a problem instance}

In this subsection, we take the optimization problem instance with $z_i=z_i(0), \forall i=1,\dots, N$ to examine the performance of Algorithm \ref{algo:DO}. For each constraint $m=1,2$, we take 
\begin{equation*}
    P^{[m]} = \begin{bmatrix}
        2/3 & 1/3 & 0 & 0 \\
        1/3 & 1/3 & 1/3 & 0\\
         0  & 1/3 & 1/3 & 1/3 \\
         0  &  0  & 1/3 & 2/3
    \end{bmatrix}
\end{equation*}
which is compatible with the line graph among the agents. 
Both Assumptions \ref{assump:network} and \ref{assump:licq} are met in this problem.
The auxiliary variables are initialized as $0$, and the parameter $\gamma$ is chosen as $0.02$. We use the quadratic programming solver in MATLAB to solve the global problem \eqref{eq:CBF_constraint} for benchmarking and also the local subproblems in Algorithm \ref{algo:DO}.

The trajectories of cost error, primal variable error,  value of coupling constraint function, and consensus error of dual variables are plotted in Figs. \ref{fig:cost}-\ref{fig:dual}, respectively. The results demonstrate the convergence of the algorithm. Particularly, we observe from Fig. \ref{fig:constraint} that the two coupling constraints are all-time satisfied, and one of them becomes active gradually. Accordingly, in Fig. \ref{fig:constraint} the local dual variables associated with the active constraint asymptotically reach consensus, and the dual variables corresponding to the inactive constraint are uniformly $0$.

\begin{figure}[htbp]
\centering
\includegraphics[width=3.3in]{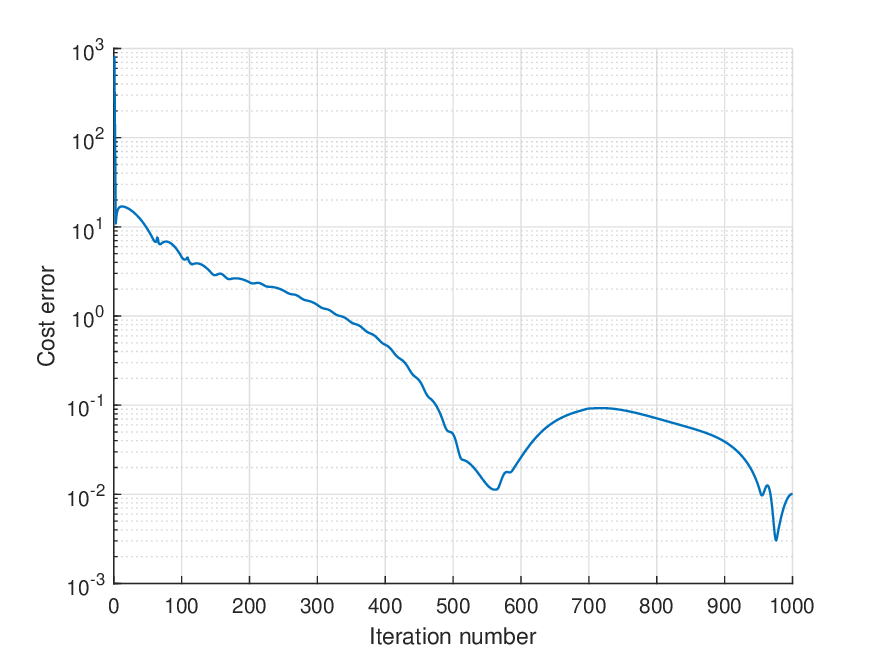}
\caption{Convergence of objective error $\phi(y)-\phi(y^*)$.}
\label{fig:cost}
\end{figure}
\begin{figure}[htbp]
\centering
\includegraphics[width=3.3in]{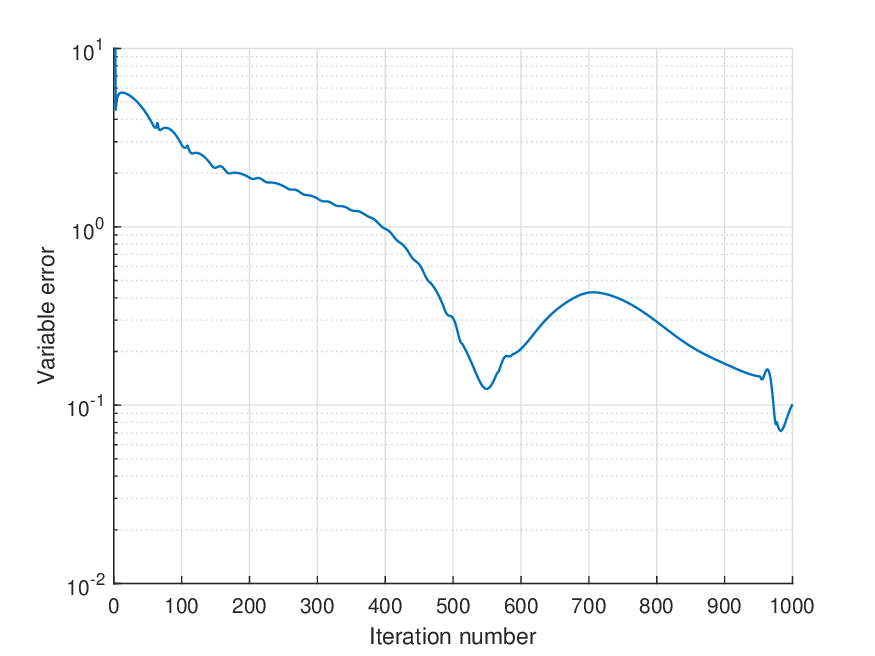}
\caption{Convergence of primal variable error $\sqrt{\sum_{i=1}^N \lVert x_i-x_i^* \rVert^2}$.}
\label{fig:variable}
\end{figure}

\begin{figure}[htbp]
\centering
\includegraphics[width=3.3in]{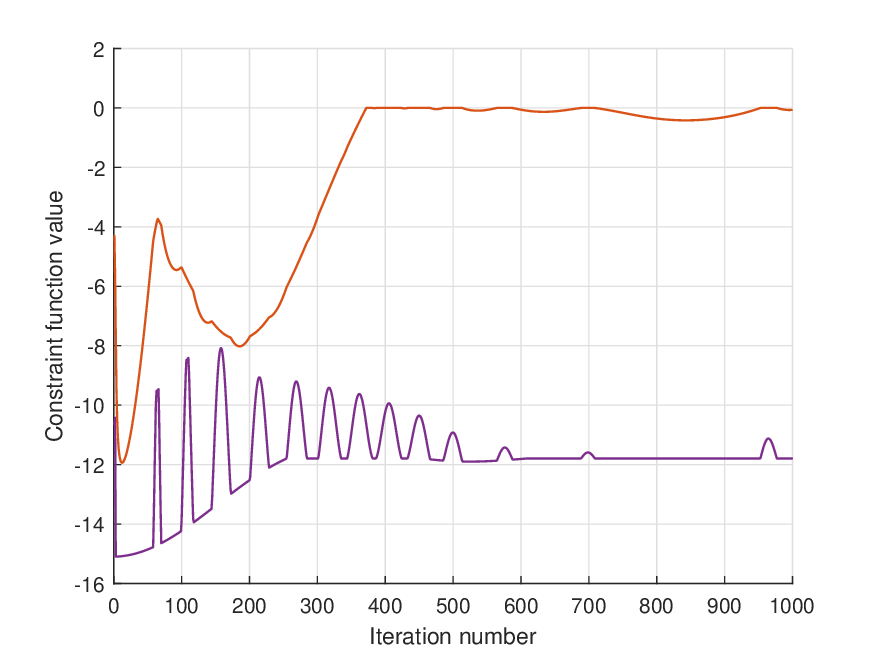}
\caption{Vector values of coupling constraint function $\sum_{i=1}^N A_i^Tx_i+b_i$.}
\label{fig:constraint}
\end{figure}
\begin{figure}[htbp]
\centering
\includegraphics[width=3.3in]{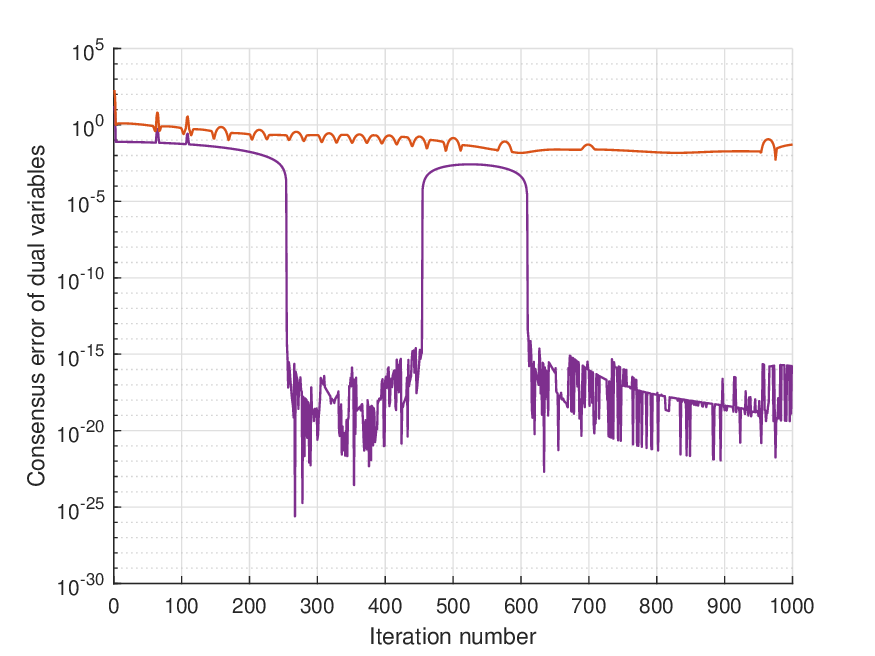}
\caption{Consensus error of local dual variables $[\lVert(I-P^{[1]}) \mu^{[1]}  \rVert, \lVert(I-P^{[2]}) \mu^{[2]} \rVert]$.}
\label{fig:dual}
\end{figure}

\subsection{Performance of the closed-loop system}

In this subsection, we contrast the performance of the CBF controller under the centralized quadratic programming solver in MATLAB and the distributed solver in Algorithm \ref{algo:DO}.

For the continuous-time system dynamics  \eqref{eq:system_dynamics}, we discretize it with a sampling period $0.01$s based on the Euler method. 
Regarding the implementation of the controller \eqref{eq:CBF_constraint} under Algorithm \ref{algo:DO}, $10$ iterations of Algorithm \ref{algo:DO} with $\gamma=0.01$ are conducted during each sampling period. In the initial iteration of each sampling period, we set the auxiliary variables to be zero. We remark that Assumption \ref{assump:network} is universally valid, and Assumption \ref{assump:licq} is only violated in exceptional cases when the coefficients associated with $x_4$ in the constraint become linearly dependent.

Due to the violation-free feature of Algorithm \ref{algo:DO}, the CBF constraints in \eqref{eq:const1} and \eqref{eq:const2} are all-time satisfied.


The comparison results are presented in Figs. \ref{fig:cent_CBF} and \ref{fig:dist_CBF}. Note that the two coupling state constraint are initially violated. The blue and red curves stand for contours of functions $(z_i^{[1]})^2 + (z_i^{[2]})^2$ and $(z_i^{[1]}-2)^2 + (z_i^{[2]}-2)^2$ with levels $1$ and $4$, respectively. Under both the centralized and distributed controllers, the agents eventually meet at one of the intersections of the two curves, suggesting that constrained consensus is achieved.


\begin{figure}[htbp]
\centering
\centering
\includegraphics[width=3.3in]{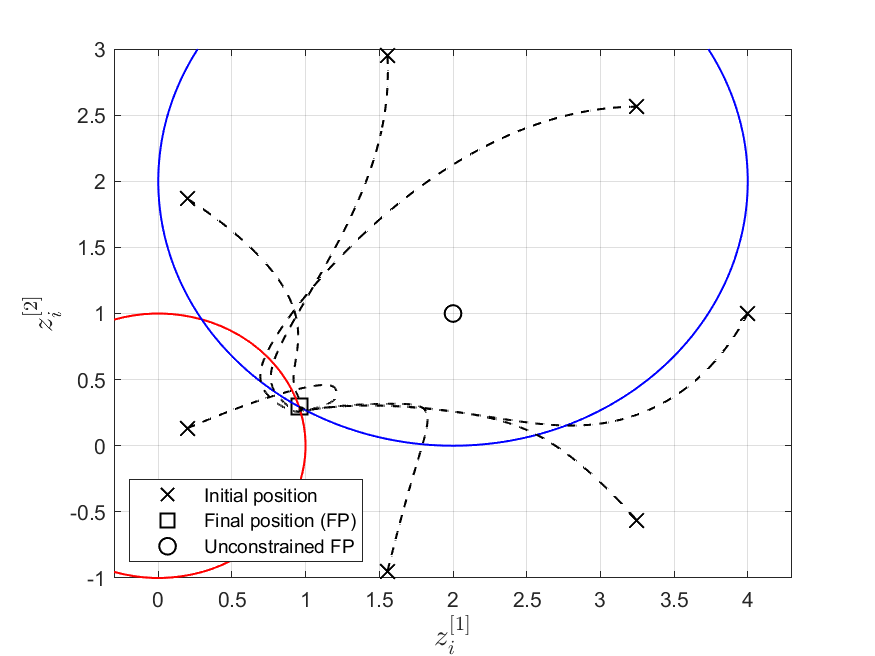}
\caption{Performance under the centralized CBF controller.}
\label{fig:cent_CBF}
\end{figure}
\begin{figure}[htbp]
\centering
\includegraphics[width=3.3in]{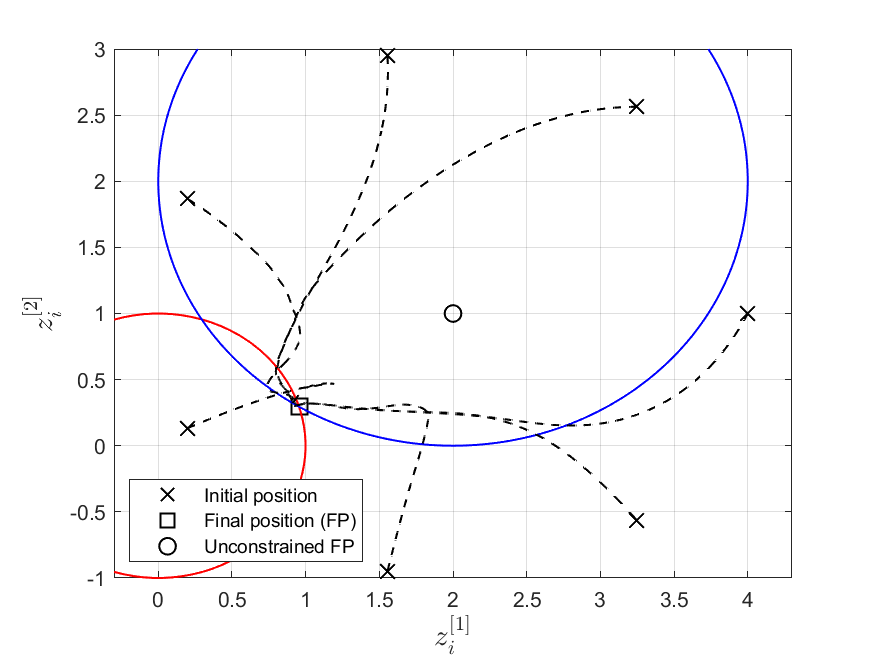}
\caption{Performance under the distributed CBF controller.}
\label{fig:dist_CBF}
\end{figure}

\section{Conclusion}

This work presented a violation-free distributed optimization algorithm for networked optimization problem with coupling constraints. Convergence of the proposed algorithm was proven under suitable conditions. We studied its application to the distributed implementation of CBF controllers where constraint satisfaction is essential, and the results demonstrated the effectiveness of the approach. 
This work paves the way for numerous opportunities for future research, including exploring the relaxation of the rank condition in Assumption \ref{assump:licq}-ii) and extending the approach to address nonlinear coupling constraints.



\bibliographystyle{plain}
\bibliography{autosam.bib}           



\end{document}